\documentclass[11pt,leqno]{article}
\begin{document}

\newtheorem{theorem}{Theorem}[section]
\newtheorem{corollary}[theorem]{Corollary}
\newtheorem{lemma}[theorem]{Lemma}
\newtheorem{proposition}[theorem]{Proposition}
\newtheorem{example}[theorem]{\sl Example}
\newtheorem{remark}[theorem]{\sl Remark}
\newtheorem{definition}[theorem]{\sl Definition}
\newcommand{\qed}{\ \ \rule{1ex}{1ex}}
\newcommand{\lc}{\left\lceil}
\newcommand{\lf}{\left\lfloor}
\newcommand{\rc}{\right\rceil}
\newcommand{\rf}{\right\rfloor}
\newcommand{\proof}{\par\noindent  {\bf  Proof.  }}
\newcommand{\bd}[1]{{\bf  #1}}
\newcommand{\RR}{{\bf  R}}
\newcommand{\NN}{{\bf  N}}
\newcommand{\ZZ}{{\bf  Z}}
\newcommand{\CC}{{\bf  C}}
\newcommand{\QQ}{{\bf  Q}}
\newcommand{\pp}{{\bf  p}}
\newcommand{\xx}{{\bf  x}}
\newcommand{\yy}{{\bf  y}}
\newcommand{\zz}{{\bf  z}}
\newcommand{\ww}{{\bf  w}}
\newcommand{\Sj}{{{\cal S}_j}}
\newcommand{\Pt}{{\tilde{P}}}
\newcommand{\Var}{{\rm Var}}
\newcommand{\diag}{{\rm diag}}
\newcommand{\Exp}{{\rm Exp}}
\newcommand{\ave}{{\rm E}}
\newcommand{\id}{{\rm id}}
\newcommand{\rev}{{\rm rev}}
\newcommand{\Xt}{{\tilde{X}}}
\newcommand{\SSS}{{S}} 
\newcommand{\gxy}{g_{xy}}
\newcommand{\Gamxy}{{\Gamma(x, y)}}
\newcommand{\vece}{{\vec{e}}}
\newcommand{\vecep}{{\vece\,'}}
\newcommand{\Qe}{{Q(\vece\,)}}
\newcommand{\Qep}{{Q(\vecep)}}
\newcommand{\we}{{w(\vece\,)}}
\newcommand{\wep}{{w(\vecep)}}
\newcommand{\Kwe}{{K_w(\vece\,)}}
\newcommand{\phie}{{\phi(\vece\,)}}
\newcommand{\psie}{{\psi(\vece\,)}}
\newcommand{\be}{{b(\vece\,)}}
\newcommand{\vecE}{{\vec{E}}}
\newcommand{\GamxyQ}{{|\Gamma(x, y)|_Q}}
\newcommand{\Gamxyw}{{|\Gamma(x, y)|_w}}
\newcommand{\gamstar}{{\gamma_*}}
\newcommand{\gamstart}{{\gamma^t_*}}
\newcommand{\gamstarone}{{\gamma^1_*}}
\newcommand{\gamstartwo}{{\gamma^2_*}}
\newcommand{\gamstarthree}{{\gamma^3_*}}
\newcommand{\gamstarfour}{{\gamma^4_*}}
\newcommand{\bte}{{b^t(\vece\,)}}
\newcommand{\bonee}{{b^1(\vece\,)}}
\newcommand{\btwoe}{{b^2(\vece\,)}}
\newcommand{\bthreee}{{b^3(\vece\,)}}
\newcommand{\bfoure}{{b^4(\vece\,)}}
\newcommand{\bthreeae}{{b^{3\mbox{\scriptsize a}}(\vece\,)}}
\newcommand{\bthreebe}{{b^{3\mbox{\scriptsize b}}(\vece\,)}}
\newcommand{\bfourae}{{b^{4\mbox{\scriptsize a}}(\vece\,)}}
\newcommand{\bfourbe}{{b^{4\mbox{\scriptsize b}}(\vece\,)}}
\newcommand{\Ac}{{\cal A}}
\newcommand{\Bc}{{\cal B}}
\newcommand{\Cc}{{\cal C}}
\newcommand{\Ec}{{\cal E}}
\newcommand{\Fc}{{\cal F}}
\newcommand{\Gc}{{\cal G}}
\newcommand{\Lc}{{\cal L}}
\newcommand{\Nc}{{\cal N}}
\newcommand{\Sc}{{\cal S}}
\newcommand{\Ft}{\Fc_t}
\newcommand{\Fi}{\Fc_{\infty}}
\newcommand{\Fn}{\Fc_n}
\newcommand{\Tt}{{\tilde{T}}}
\newcommand{\ft}{{\tilde{f}}}
\newcommand{\Jt}{{\tilde{J}}}
\newcommand{\tp}{t+}
\newcommand{\tm}{t-}
\newcommand{\vdn}{\| \pi_n - \pi \|}
\newcommand{\vdt}{\| \pi_t - \pi \|}
\newcommand{\sep}{{\rm sep}}
\newcommand{\Fets}{\Fc^*_{=t}}
\newcommand{\ppi}{{\mbox{\boldmath $\pi$}}}
\newcommand{\PP}{{\bf  P}}
\newcommand{\PPh}{\PP^{(h)}}
\newcommand{\GG}{{\bf  G}}
\newcommand{\II}{{\bf  I}}
\newcommand{\XX}{{\bf  X}}
\newcommand{\Qsp}{{Q^*}'}
\newcommand{\sn}{\sigma_n}
\newcommand{\tn}{\tau_n}
\newcommand{\tnm}{\tau_{n-1}}
\newcommand{\gt}{{\tilde{g}}}
\newcommand{\Gt}{{\tilde{G}}}
\newcommand{\pimin}{\pi_{\min}}
\newcommand{\bege}{\begin{equation}}
\newcommand{\ene}{\end{equation}}
\newcommand{\begp}{\begin{proposition}}
\newcommand{\enp}{\end{proposition}}
\newcommand{\begt}{\begin{theorem}}
\newcommand{\ent}{\end{theorem}}
\newcommand{\begl}{\begin{lemma}}
\newcommand{\enl}{\end{lemma}}
\newcommand{\begc}{\begin{corollary}}
\newcommand{\enc}{\end{corollary}}
\newcommand{\begr}{\begin{remark}}
\newcommand{\enr}{\end{remark}}
\newcommand{\begd}{\begin{definition}}
\newcommand{\enf}{\end{definition}}
\newcommand{\begx}{\begin{example}}
\newcommand{\enx}{\end{example}}
\newcommand{\bega}{\begin{array}}
\newcommand{\ena}{\end{array}}
\newcommand{\Line}{$\underline{\qquad\qquad\qquad\qquad}$}
\newcommand{\X}{$\qquad$}
\newcommand{\lb}{\left\{}
\newcommand{\rb}{\right\}}
\newcommand{\lsb}{\left[}
\newcommand{\rsb}{\right]}
\newcommand{\lp}{\left(}
\newcommand{\rp}{\right)}
\newcommand{\ls}{\left|}
\newcommand{\rs}{\right|}
\newcommand{\lss}{\left\|}
\newcommand{\rss}{\right\|}
\newcommand{\vs}{\vspace{\smallskipamount}\noindent}
\newcommand{\vm}{\vspace{\medskipamount}\noindent}
\newcommand{\vb}{\vspace{\bigskipamount}\noindent}
\newcommand{\ra}{\rightarrow}
\newcommand{\la}{\leftarrow}
\newcommand{\implies}{\Longrightarrow}
%


\newpage

\section*{An interesting spectral gap problem, from Jim Fill \\ 
(more formally, from \\ James Allen Fill, \\ The Johns Hopkins University, \\ Department of Applied Mathematics and Statistics)} \label{s:gap}

\subsection*{Preamble} \label{ss:preamble}

\hspace\parindent
At the request of Laszlo Babai, founder and an editor of the free online journal
``Theory of Computing'' (ToC), {\tt theoryofcomputing.org},\ in August, 2025, I am posting on the arXiv,
essentially unedited and not updated, a combination of two closely related sets of unpublished notes from 2003.
ToC is keen on publishing links to all bibliography items, and a paper soon to be published there makes progress
on a conjecture in my 2003 notes.
The sections ``The problem'', ``Evidence in favor of the conjecture'', ``Facts about the spectral structure of the matrix~$K$'', and ``Stronger conjectures'' previously formed a document entitled ``An interesting spectral gap problem, from Jim Fill''.
The sections ``Introduction:\ Self-organizing lists'' and ``The move-ahead-$1$ (MA1) rule'' formed a document entitled ``Background on the gap problem''.  The two sets of notes have inspired some research, including (as one example, with no attempt here at a literature survey)~\cite{BMRS}.

\subsection*{The problem} \label{ss:problem}

\hspace\parindent
Here is an interesting, simple-to-state problem about spectral gaps that I have no idea
how to solve.  I would love to receive help---or a complete solution!  I will provide my
motivation for the problem to anyone who asks.

Let~$n \geq 2$ be an integer, and for $i, j \in \{1, \ldots, n\}$ with $i \neq
j$, let $p_{i j} \in (0, 1)$.  We assume that $p_{i j} + p_{j i} = 1$ for all $i, j$. 
Given these parameters, we consider the following informal description of a specific
Markov chain:
\begin{quote}
Choose a position $r \in \{2, \ldots, n\}$ uniformly at random; with probability
$p_{x_r, x_{r - 1}}$ move from the current permutation $x = (x_1, \ldots, x_n)$ to
$$
(x_1, \ldots, x_{r-2}, x_r, x_{r-1},x_{r+1}, \ldots, x_n),
$$
and with the remaining probability stay at~$x$.
\end{quote}
To describe the transition matrix~$K$ of this chain more formally,
\begin{itemize}
\item if~$x$ and~$y$ differ precisely by a transposition of adjacent items~$i$ and~$j$,
with~$i$ to the (immediate) left of~$j$ in~$x$ and~$j$ to the (immediate) left of~$i$
in~$y$, then $K_{x y} = p_{j i} / (n - 1)$;
\item if not, and if $x \neq y$, then $K_{x y} = 0$;
\item the diagonal values are $K_{x x} = 1 - \sum_{y \neq x} K_{x y}$.
\end{itemize} 

As examples, if $n = 2$ and the rows and columns are given in the order $(1, 2), (2, 1)$,
then
$$
K = \left[ \begin{array}{cc}
                       p_{1 2} & p_{2 1} \\
                       p_{1 2} & p_{2 1}
                      \end{array}
               \right];
$$
and if $n = 3$ and the rows and columns are given in the order $(1, 2, 3), (2, 1, 3),
(1, 3, 2)$, $(2, 3, 1), (3, 1, 2), (3, 2, 1)$, then
$$
K = \frac{1}{2} \left[ \begin{array}{cccccc}
        p_{1 2} + p_{2 3} & p_{2 1} & p_{3 2} & 0 & 0 & 0 \\
        p_{1 2} & p_{2 1} + p_{1 3} & 0 & p_{3 1} & 0 & 0 \\
        p_{2 3} & 0 & p_{1 3} + p_{3 2} & 0 & p_{3 1} & 0 \\
        0 & p_{1 3} & 0 & p_{2 3} + p_{3 1} & 0 & p_{3 2} \\
        0 & 0 & p_{1 3} & 0 & p_{3 1} + p_{1 2} & p_{2 1} \\
        0 & 0 & 0 & p_{2 3} & p_{1 2} & p_{3 2} + p_{2 1}
                      \end{array}
               \right].
$$

It is easy to see that~$K$ is reversible with respect to its stationary distribution
$$
\pi_x = z^{-1} \prod_{r, s:\ 1 \leq r < s \leq n} p_{x_r, x_s},
$$
where~$z$ is a normalizing constant chosen so that $\sum_x \pi_x = 1$; that is,
$$
\pi_x K_{x y} = \pi_y  K_{y x}\mbox{\ \ for all $x, y$.}
$$
Clearly, all the eigenvalues of~$K$ lie in $(-1, 1]$, and~$1$ is a simple eigenvalue.  (For
a strengthening, see Theorem~2 below.)
\smallskip

\par\noindent
{\bf Definition.}\ Call the parameter vector $\pp = (p_{i j}: 1 \leq i < j \leq n)$
\emph{regular} if
\begin{enumerate}
\item $p_{i - 1, i} \geq 1 / 2$ for $2 \leq i \leq n$,
\item $p_{i - 1, j} \geq p_{i j}$ for $2 \leq i < j \leq n$, and 
\item $p_{i, j + 1} \geq p_{i j}$ for $1 \leq i < j \leq n - 1$.
\end{enumerate}
\smallskip

Define the (spectral) \emph{gap}~$\lambda_K$ for~$K$ by $\lambda_K := 1 - \beta_K$, where
$\beta_K$ is the second-largest eigenvalue of~$K$.
\bigskip

\par\noindent
{\bf Gap Problem.}\ Prove (or disprove) the following conjecture: Among all regular choices
of~$\pp$, the gap~$\lambda_K$ is minimized by choosing $p_{i j} = 1 / 2$ for all $1 \leq i
< j \leq n$.
\medskip

In lieu of a solution to the Gap Problem, I would be nearly as happy with any upper bound on
the set
$$
R_n := \{ 1 / \lambda_K: \mbox{$\pp$ is regular} \}
$$
that grows only polynomially in~$n$.  (Cf.\ Theorem~1 below.)  

\subsection*{Evidence in favor of the conjecture} \label{ss:evidence}

\hspace\parindent
(a) The conjecture is true for $n = 2$, for which $\lambda_K = 1$ for \emph{any} choice of
$p_{1 2}$.
\medskip

(b) The conjecture is true for $n = 3$.  I used {\tt Mathematica} to discover the
eigenstructure of~$K$ in this case, and the results are easily verified.  In particuar, we
have the formula
$$
\lambda_K = \frac{1}{2} \left[ 1 - \left( p_{1 2} p_{2 3} p_{3 1} + p_{3 2} p_{2 1} p_{1
3} \right)^{1 / 2} \right],
$$
and $\beta_K$ has multiplicity~$2$ as an eigenvalue.  As~$\pp$ ranges over the regular
vectors, the gap~$\lambda_K$ is nondecreasing in each of $p_{1 2}$, $p_{2 3}$, and $p_{1 3}$
as the other two $p$'s are held fixed.  In particular, one finds that $\lambda_K$ is
minimized by a given choice~$\pp$ if and only if $p_{1 2} = 1 / 2 = p_{2 3}$ (with the
value of $p_{1 3}$ then being curiously irrelevant).
\medskip

(c) There are strong numerical indications that the conjecture is true for $n = 4$ and for
$n = 5$.  For $n = 4$ (respectively, $n = 5$), I calculated $\lambda_K$ (numerically) for
all regular~$\pp$ for which each $p_{i j}$ is a multiple of $0.05$ (respectively, $0.10$),
and the minimum value occurred, as conjectured, in the unweighted case $p_{i j} \equiv
0.50$. 

\subsection*{Facts about the spectral structure of the matrix~$K$} \label{ss:theorems}

\hspace\parindent
In this section I'll present a few facts about the spectral structure of the matrix~$K$. 
The first is a result from~\cite{lozenge}.
\medskip

\par\noindent
{\bf Theorem~1.}\ In the unweighted case $p_{ij} \equiv 1 / 2$ the gap is
$$
\lambda_K = (1 - \cos(\pi / n)) / (n - 1).
$$

Notice that $(1 - \cos(\pi / n)) / (n - 1) \sim \pi^2 / (2 n^3)$ as $n \to \infty$.

In the unweighted case, the multiplicity of the eigenvalue $\beta_K = 1 - \lambda_K$ is at
least (in fact, I believe, is exactly) $n - 1$.  Indeed, for each $1 \leq i \leq n$, the
label~$i$ performs a (slowed-down) random walk on the path $\{1, \ldots, n\}$, with
spectral gap equal to $(1 - \cos(\pi / n)) / (n - 1)$.  Further, the corresponding
walk-eigenvector lifts to an eigenvector with eigenvalue $\beta_K$ for the permutation
chain.  These~$n$ eigenvectors have vanishing sum, but any $n - 1$ of them are linearly
independent.

Next I present a fact about the spectral structure of~$K$ for general (not necessarily
regular) $\pp$ that I first noticed from numerical computations.
\medskip

\par\noindent
{\bf Theorem~2.}\ The matrices~$K$ and $I - K$ are similar.  Therefore, the eigevalues
of~$K$ occur in pairs summing to unity, and the spectrum of~$K$ lies in $[0, 1]$.
\medskip

\par\noindent
\emph{Proof.}\ List the permutations of length~$n$ in an order so that, for every~$k$, the
$k$th-to-last permutation is the reverse of the $k$th permutation, as I did above in
writing out~$K$ for $n = 3$.  Let $D := \diag(\pi)$ (recalling that~$\pi$ is the stationary
distribution for~$K$), let~$B$ denote the backward identity matrix (the matrix with $b_{i,
n + 1 - i} = 1$ for all~$i$ and $b_{i j} = 0$ otherwise), and let~$E :=
\diag(\mbox{sgn}(x))$, where $\mbox{sgn}(x)$ denotes the sign of the permutation~$x$.  Set
$C: = E B D$.  Then~$C$ is invertible and
$$
I - K = C K C^{-1}.~\qed
$$

\subsection*{Stronger conjectures} \label{ss:stronger}

\hspace\parindent
In this section I'll discuss natural stronger conjectures than the one proposed in the Gap
Problem.
\medskip

(a) Based on evidence~(b) on page~2 above, it was natural to guess for general~$n$ that,
as~$\pp$ ranges over the regular vectors, the gap~$\lambda_K$ is nondecreasing in each
$p_{ij}$ as the other parameters~$p_{rs}$ are held fixed.  But this is false for $n = 4$. 
For example, when $p_{1 2} = 0.5$, $p_{2 3} = 0.7$, $p_{3 4} = 0.5$, $p_{1 3} = 0.7$,
$p_{2 4} = 0.8$, and $p_{1 4} = 0.9$, we have $\lambda_K \approx 0.1261$, but when $p_{1 2}$
is increased to $0.6$ we have the smaller gap $\lambda_K \approx 0.1259$.

But for evidence of a different sort of monotonicity conjecture, see~(d) below.
\medskip

(b) It is also false (even for $n = 4$) that every eigenvalue in the upper half of the
spectrum of~$K$ is maximized in the unweighted case.  Indeed, the fifth-largest eigenvalue
is approximately $0.7887$ in the unweighted case but is about $0.7944$ when
$p_{1 2} = p_{2 3} = p_{3 4} = p_{1 3} = 0.50$ and $p_{2 4} = p_{1 4} = 0.95$.
\medskip

(c) Numerical evidence suggests that, for regular~$\pp$, the gap~$\lambda_K$ achieves the
minimum value of Theorem~1 if and only if there exists $1 \leq i \leq n$ such that
$p_{i j} = 1 / 2$ for all $j \neq i$.  Further, in such cases the multiplicity (call
it~$\mu$) of the eigenvalue~$\beta_K$ appears to equal the number (call it~$\nu$) of such
values~$i$, with two exceptions: (i) if $\nu = n$ (i.e.,\ in the unweighted case), then
$\mu = n - 1$; and (ii) if $\nu = n - 2$ (i.e.,\ if $p_{1 n} = 1 - p_{n 1} > 1 / 2$ and
$p_{i j} = 1 / 2$ otherwise), then $\mu = n - 1$.
\medskip

(d) Observe that the regular vectors form a convex set.  Let~$\pp_*$ denote~$\pp$ in the
unweighted case.  For given regular~$\pp$ and $0 \leq t \leq 1$, define $\pp(t) := (1 - t)
\pp + t \pp_*$, and let $K(t)$ denote the corresponding transition matrix.  The Gap Problem
formulates the conjecture that $K(1)$ has a smaller gap than does $K(0)$.  Numerical
evidence leads to the stronger conjecture that the gap $\lambda_{K(t)}$ is nonincreasing
(and also convex) in $t \in [0, 1]$.

\subsection*{Introduction: Self-organizing lists} \label{ss:lists}

\hspace\parindent
Here is background about why I am interested in the gap problem.  In brief, the problem
arises in connection with the analysis of certain Markov chains on self-organizing
lists.  (I would be happy to elaborate on this background for interested parties.)

Suppose that we have a set of $n$ {\em records\/} (e.g.,\ items of data in a computer),
each equipped with an identifying {\em key\/}, or label, and suppose that we wish to
store these records in a way that facilitates searching and sorting.  One frequently
employed data structure for this is a list that is ordered linearly, say from ``front''
to ``back.''  When a record is requested, the list is searched sequentially from front
to back until the desired key is located.  (For simplicity, we shall assume that only
records appearing in the list are ever requested.)  The cost of meeting a request can be
defined as the number of records that must be examined to reach the requested record, or
more generally as some function of this number.  Such a list is called a {\em linear
search list\/}.  Any standard text on data structures, e.g.,~\cite{Knuth}, \cite{LD},
\cite{Mehlhorn}, or~\cite{van Leeuwen}, may be consulted for a general introduction.

Without loss of generality, we write $[n] := \{1, 2, \ldots, n\}$ for the set of keys. 
Suppose now that independently at each unit of time a record is requested (accessed)
according to some fixed probability distribution; let the ``weight'' $w_i$ denote the
probability that item $i$ is requested.  It is clear that the optimal static order
in which to arrange the items is from front to back in decreasing (by which we shall
mean nonincreasing) order of the $w_i$'s.  However, typically the long-run
usage frequencies $w_i$ are unknown.  One remedy that immediately leaps to mind
is to keep a running counter of the actual frequency of access for each record and
dynamically reorder the items in decreasing order of frequencies.  By the strong law of
large numbers, the order obtained by this counter scheme will converge with
probability~1 to the optimal static order, as observed in~\cite{Rivest}.  However, this
analysis fails to take into account the costs of the rearrangements and the memory
requirements of the counters.  

A commonly used alternative is to ``let the records organize themselves.''  The
idea of a {\em self-organizing (linear search) list\/} (SOL) scheme is to rearrange the
records after each access, according to some well-specified rule without the need
for counters, in such a way that frequently accessed records tend to gravitate toward
the front of the list, thereby reducing search costs.  If the (conditional)
probabilities of the various possible rearrangements at a given epoch depend solely on
the order of the items at that epoch and the item selected, then the time-dependent
order of the items forms a Markov chain on the symmetric group $S_n$.  A major goal
is to assess rates of convergence to stationarity and other operating characteristics
for the permutation chain.  

By far, the two most frequently used and most frequently discussed schemes are the
simple {\em move-to-front\/} (MTF) and {\em transposition\/}, or {\em
move-ahead-$1$\/} (MA1), rules.  My gap problem evolved from consideration of the MA1
rule, so I will almost exclusively consider that rule here.   Consider a fixed epoch
at which the records are in order $x = (x_1, \ldots, x_n)$, with $x_j$ denoting
the key (label) of the record in position $j$.  Under the MTF rule, if the record in
position $i$ is requested (this happens with probability $w_{x_i}$), it is brought to
the front and the order is subsequently rearranged to $(x_i, x_1, \ldots, x_{i-1},
x_{i+1}, \ldots, x_n)$.  Using MA1, the selected item is instead moved one position
ahead, resulting in $(x_1, \ldots, x_{i-2}, x_i, x_{i-1},x_{i+1}, \ldots,
x_n)$; if $x_1$ is selected, the order remains at $x$.

A gentle introduction to self-organizing Markov chains is provided
by~\cite{Hendricks89}.  Hester and Hirschberg~\cite{HH} survey the state of the field,
including known results and open problems, as of~1985.  We are concerned only with {\em
probabilistic\/} analysis of self-organizing (sequential) search (SOS).  Alternatively,
one may employ {\em amortized analysis\/} (see~\cite{ST1} and~\cite{ST2}) which produces
a sort of ``worst case average performance'' of a self-organizing scheme.

{\em Note\/}:  Without any real loss of generality we suppose for the remainder of this
document that $w_1 \geq w_2 \geq \cdots \geq w_n > 0$.

\subsection*{The move-ahead-$1$ (MA1) rule} \label{ss:tr}

\hspace\parindent
The MA1 rule moves records far less than does MTF.  Intuitively, therefore, MA1
performs better than MTF when the records are in ``good'' order and not as well when
the order is ``bad.''  (One also expects that MA1 will not perform as well as MTF in the
Markov requests model when there is positive dependence among successive requests.) 
Rivest~\cite{Rivest} showed that {\em stationary\/} expected search cost (ESC) is smaller
for MA1 than for MTF, but little has been done analytically to show that MA1 peforms
poorly due to slow convergence when the initial order is bad.  The summary of the
literature on MA1 by Hester and Hirschberg~\cite{HH} is accurate:
\begin{quote}
{\em Move-to-front\/} is the only permutation algorithm that has been extensively
analyzed.  {\em Transpose\/} is widely mentioned in the literature, but authors merely
state an inability to obtain theoretical bounds on its performance, and most of the
work on other algorithms merely shows how their behavior changes as parameters vary. 
Direct bounds on the behaviors of these algorithms are needed to allow realistic
comparisons.
\end{quote}

The MA1 rule gives rise to a reversible Markov chain on $S_n$. 
Hendricks~\cite{Hendricks76} derived the stationary distribution~$\pi$ of the chain:
\begin{equation}
\label{stationary}
\pi(x) = z^{-1} w^{n - 1}_{x_1} \cdots w^0_{x_n} = z^{-1} \prod_{i = 1}^n w^{n -
i}_{x_i},
\end{equation}
where~$z$ is a normalizing constant chosen so that $\sum_x \pi(x) = 1$.
(The meaning of ``reversibility'' is that $\pi(x) M(x, y) \equiv \pi(y) M(y, x)$,
where~$M$ is the transition matrix for the chain.  Thus~$\pi$ is the unique
stationary distribution for the chain, and~$M$ is diagonally similar to the symmetric
matrix~$S$, where $S(x, y) := \pi(x)^{1/2} M(x, y) \pi(y)^{-1/2}$.  Also, $M(x, x) > 0$
for all~$x$, and so the eigenvalues of~$S$ are all in $(-1, 1]$.)  But the consequent
expression for the stationary ESC is computationally intractable.  To our knowledge, the
only estimates of the rate of convergence to stationarity for the permutation-chain that
have been obtained are for the unweighted case ($w_i \equiv 1$).  Using a clever
coupling argument for the upper bound and a ``follow-a-single-card'' approach for the
lower, Aldous~\cite{Aldous83} showed that the number of moves required for small
variation distance is between the orders of $n^3$ and
$n^3 \log n$.  Diaconis and Saloff-Coste~\cite{DSCgroup} obtained the same result by
comparing the eigenvalues of the MA1 chain with those of the ``random transpositions''
shuffle of Diaconis and Shahshahani~\cite{DSh81}.  They conjectured that the correct
order is $n^3 \log n$.  [This was eventually established by Wilson~\cite{lozenge}, who
showed that the spectral gap (difference between~$1$ and the second largest eigenvalue
of~$K$) is exactly $(1 - \cos(\pi / n)) / (n - 1)$.]  From~\cite{DFP}, the corresponding
order for MTF is only $n \log n$; thus, for equal weights, MA1 is indeed a good deal
slower to converge to uniform order.

I don't know how to extend these techniques to the weighted case.  To mention one
line of attack which has been unproductive (namely, the eigenvalue comparison technique
of Diaconis and Saloff-Coste~\cite{DSCrev}), comparison with an unweighted shuffle such
as random transpositions is ineffective because of the existence of large values of the
ratio of the two stationary distributions. 

The MA1 chain sometimes is provably slow to mix.  For example, consider the geometric
weights $w_i \equiv 2^{-i} / (1 - 2^{-n})$, and suppose that the MA1 chain is started in
a permutation with record~$n$ in front of record $n - 1$.  In stationarity, the
probability that $n - 1$ is in front of~$n$ is larger than $1/2$ (and, by the way,
strictly so for $n \geq 3$).  Yet $n - 1$ cannot possibly move ahead of record~$n$
until it is requested at least once, and the request probability $w_{n - 1}$ is
exponentially small in~$n$.

But I am willing to employ \emph{any} chain with the stationary
distribution~(\ref{stationary}), since (among other goals) I am interested in sampling
from~(\ref{stationary}) for various choices of the weights~$w$.  So I looked for a chain
with stationary distribution~(\ref{stationary}) that might be rapidly mixing regardless
of the weights.  It was my intuition that the following chain would work well. 
Choose a position $i \in \{2, \ldots, n\}$ \emph{uniformly} at random; with probability
$w_{x_i} / (w_{x_{i - 1}} + w_{x_i})$ move from the current permutation $x = (x_1,
\ldots, x_n)$ to $(x_1, \ldots, x_{i-2}, x_i, x_{i-1},x_{i+1}, \ldots, x_n)$, and with
the remaining probability stay at~$x$.  It is easy to check that this chain is
reversible with respect to the distribution~(\ref{stationary}).
More generally, suppose we introduce parameters
$p_{i,j}$ as in my gap problem and consider the following chain:
\begin{quote}
Choose a position $r \in \{2, \ldots, n\}$ uniformly at random; with probability
$p_{x_r, x_{r - 1}}$ move from the current permutation $x = (x_1, \ldots, x_n)$ to
$$
(x_1, \ldots, x_{r-2}, x_r, x_{r-1},x_{r+1}, \ldots, x_n),
$$
and with the remaining probability stay at~$x$.
\end{quote}
Then it is easy to see that this chain is reversible with respect to the distribution
\begin{equation}
\label{statgen}
\pi(x) = z^{-1} \prod_{r, s:\ 1 \leq r < s \leq n} p_{x_r, x_s},
\end{equation}
where again~$z$ is a normalizing constant.  The generalization has come by extending
$w_r / (w_r + w_s)$ to $p_{r, s}$.  Observe that the assumed decreasingness of the
weights $w_r$ implies, for that special case, the regularity property of the parameters
$p_{r, s}$ I required in my write-up of the gap problem.

If I can resolve my gap problem positively, it will have several positive
consequences.  In addition to providing a way to do approximate sampling
from~(\ref{stationary}) [or, more generally, from~(\ref{statgen})] \emph{with guaranteed
bounds on the error} and rapidly (at least when started from the identity permutation),
it will also allow me to get rapid mixing results for the MA1 chain itself, for some
choices of~$w$, by comparison methods.

Ideally, the positive resolution would allow for generalization.  Indeed, similar
mixing-rate challenges await for natural weighted versions of other classical Markov
chains, for example, the Bernoulli--Laplace model for diffusion (see~\cite{Diacbook},
page~56).  In the unweighted case, Diaconis and Shahshahani~\cite{DSh87} use
representation theory for Gelfand pairs to give precise formulas for the eigenvalues and
eigenvectors of the transition matrix and obtain sharp estimates for convergence.  In
the weighted case, only the stationary distribution is known (Persi Diaconis, personal
communication).  

\subsection*{} \label{s:bib}

\end{document}